\documentclass[11pt,oneside]{amsart}

\usepackage{amsmath,ifthen, amsfonts, amssymb}
\newcommand{\showcomments}{no}

\newsavebox{\commentbox}
\newenvironment{com}%
{\ifthenelse{\equal{\showcomments}{yes}}%
{\footnotemark
     \begin{lrbox}{\commentbox}
     \begin{minipage}[t]{1.25in}\raggedright\sffamily\tiny
     \footnotemark[\arabic{footnote}]}
{\begin{lrbox}{\commentbox}}}%
{\ifthenelse{\equal{\showcomments}{yes}}%
{\end{minipage}\end{lrbox}\marginpar{\usebox{\commentbox}}}
{\end{lrbox}}}

\newtheorem{thm}{Theorem}[section]
\newtheorem{lem}[thm]{Lemma}
\newtheorem{cor}[thm]{Corollary}

\theoremstyle{definition}

\newtheorem{prob}[thm]{Problem}

\DeclareMathOperator{\Out}{Out}

\DeclareMathOperator{\dimension}{dim}

\newcommand{\field}[1]{\mathbb{#1}}
\newcommand{\integers}{\ensuremath{\field{Z}}}

\newcommand{\reals}{\ensuremath{\field{R}}}

\newcommand{\semidirect}{\ltimes}

\begin{document}

\title[The Tits alternative for CAT(0) cubical complexes]
{The Tits alternative for CAT(0) cubical complexes}

\author[M.~Sageev]{Michah Sageev}
       \address{Dept. of Math.\\
                Technion \\
                Haifa 32000, Israel}
       \email{sageevm@techunix.technion.ac.il}
\thanks{Sageev's research supported by BSF and ISF}

\author[D.~T.~Wise]{Daniel T. Wise}
       \address{Math. \& Stats.\\
                McGill University \\
                Montreal, Quebec, Canada H3A 2K6 }
       \email{wise@math.mcgill.ca}
\thanks{Wise's research supported by  FCAR and NSERC}

\subjclass[2000]
{20F67, 
20E08 
}

\keywords{CAT(0) cubical complex, Tits alternative}
\date{\today}

\begin{com}
{\bf \normalsize COMMENTS\\}
ARE\\
SHOWING!\\
\end{com}

\begin{abstract}
We prove a Tits alternative theorem for groups acting on CAT(0) cubical
complexes. Namely, suppose  that  $G$
is a group for which there is a bound on the orders of its finite
subgroups.
We prove that if $G$ acts properly on a finite-dimensional CAT(0)
cubical complex, then either $G$ contains a free subgroup of rank~2 or
$G$ is  finitely generated and virtually abelian. In particular the
above conclusion holds for any group $G$ with a free action on a
finite-dimensional CAT(0) cubical complex.
\end{abstract}

\maketitle

\section{Introduction}

A group $G$ satisfies the {\em Tits alternative} if
for every subgroup $H$ of $G$, either $H$ is virtually solvable or
$H$ contains a  free subgroup of rank~$2$.
The Tits alternative is named for Jacque~Tits,
who discovered
that it holds for any finitely generated  linear group \cite{Tits72}.
  The Tits alternative has
  many important consequences - for instance, it implies
von Neumann's dichotomy
  between amenability and containing a  free subgroup of rank~$2$.

The purpose of this note is to prove:

\begin{thm}\label{thm:main}
Suppose that $G$ is a group for which there is a bound on the order of its  finite subgroups.
Suppose that $G$ acts properly on a CAT(0) cubical complex.  Then for
each subgroup $H$ of $G$, either $H$ contains a rank~$2$ free subgroup,
or $H$ is virtually a finitely generated abelian group.
\end{thm}

Since finite subgroups of $G$ have fixed points in any action of
$G$ on a CAT(0) space, we have the following corollary.

\begin{cor}\label{cor:main}
Let $G$ be a group which acts freely on a finite dimensional CAT(0)
cubical complex.
Then for each subgroup $H$ of $G$, either
$H$ contains a rank~$2$ free subgroup,
or $H$ is virtually a finitely generated abelian group.
\end{cor}

Theorem \ref{thm:main} is false without the finite dimensional hypothesis.
  Indeed, Thompson's group is not virtually solvable
  and contains no rank~$2$ free subgroup \cite{BrinSquier85},
but Thompson's group acts freely on an infinite
dimensional (yet proper) CAT(0) cubical complex \cite{Farley2003}.

Theorem~\ref{thm:main} also fails to hold without a bound on the
order of finite subgroups. Indeed, let $G$ be the ascending union
of finite groups $G=G_1\subset G_2\subset \cdots$ (e.g. the direct
product of countably many finite groups). Then there exists a
``coset tree" $T$ on which $G$ acts (see \cite{Serre80}); namely,
the lattice of left cosets of the subgroups $G_i$ is a tree. The
natural left action of $G$ on $T$ is proper, but $G$ does not
satisfy the conclusion of Theorem \ref{thm:main}.

We list some cases in which the Tits alternative is already known
for groups arising in geometric group theory. Note that $G$
satisfies the {\em weak Tits alternative} if either $G$ is
virtually solvable, or $G$ contains $F_2$ (as opposed to requiring
this for each finitely generated subgroup).

\medskip

\underline{The Tits Alternative} holds for the following groups:
\begin{enumerate}
\item word-hyperbolic groups \cite{Gromov87}
\begin{com} I think its there...\end{com}
\item
$\Out(F_n)$ \cite{BestvinaFeighnHandel2000}
\item
Foldable cubical chamber complexes \cite{BallmannSwiatkowski99}
\item
$\pi_1M$ where $M$ is a nonpositively curved real-analytic $4$-manifold
satisfying a certain cycle condition
\cite{Xiangdong2004}
\item
fundamental groups of CAT(0) square complexes with no fake planes
\cite{Xiangdong05}
\begin{com}
And I proved it for $\pi_1X$ if $X$ has nonpositive immersions.
\end{com}

\medskip
\noindent
\hspace{-.9cm}\underline{The Weak Tits Alternative} holds
for $\pi_1$ of a compact
space satisfying:

\medskip
\item
$C(4)$-$T(4)$: \cite{Collins73}
\item
$C(6)$:  \cite{ElMosalamy86}
(further discusses the $C(4)$-$T(4)$ case
and references $C(6)$ case to dissertation
in bibliography)
\item
$C(3)$-$T(6)$: \cite{EdjvetHowie88}
\item
Finite 2-complexes satisfying the Gersten-Pride weight
test with positive rational angles \cite{CorsonTrace96}
\begin{com}
I need to check if this is not an overstatement.
\end{com}
\item
Finite nonpositively curved $2$-complexes \cite{BallmannBrin95}
\end{enumerate}

The results most closely related to this paper are
those of Ballman and \'{S}wi\c{a}tkowski  \cite{BallmannSwiatkowski99} and Xie \cite{Xiangdong05}.
We note that Ballmann and \'{S}wi\c{a}tkowski actually proved
a stronger result allowing larger stabilizers in a more restrictive setting.
We also note that Theorem~\ref{thm:main} unifies some of the results
described above in the small-cancellation case.
Indeed, fundamental groups of compact $C'(\frac16)$ and
$C'(\frac14)$-$T(4)$ complexes
act with uniformly bounded stabilizers on
finite dimensional CAT(0) cube complexes \cite{WiseSmallCanCube}.
These are the most important metric subcases
of the $C(6)$ and $C(4)$-$T(4)$ groups above.
\begin{com}
Increase my citation impact factor!
\end{com}

We close the introduction with the following problems:
\begin{prob}
Does the Tits alternative hold for $G$ if:

\begin{enumerate}
\item  $G$ acts properly discontinuously and
cocompactly on a CAT(0) space? \\
(This is even unknown for CAT(0) manifolds.)

\item$G$ is automatic?
\end{enumerate}
\end{prob}

\section{Preliminaries}
\subsection{Cubical complexes and hyperplanes}

We recall some basic facts about hyperplanes (see \cite{Sageev95}
for details). Let $X$ be a CAT(0) cubical complex. Two
edges $e$ and $f$ are {\em square equivalent} if they
are opposite edges of a
square $\sigma$ in $X$.
  We let $\sim$ denote the equivalence relation generated
by square equivalence and let $\overline{e}$ denote the
equivalence class containing $e$. We identify each $n$-cube $\sigma$
of $X$ with the standard unit cube in $\reals^n$.
  A {\em midcube} of $\sigma$ is the intersection of $\sigma$
with an $(n-1)$-dimensional hyperplane parallel to one of the faces
of $\sigma$ and containing the center of $\sigma$. A {\it
hyperplane} in $X$ is the union of all midcubes meeting a
particular equivalence class of edges. We note some basic facts
about hyperplanes.

\begin{enumerate}
     \item A hyperplane
     meets each cube in at most one midcube.
     \item A hyperplane separates $X$ into
     precisely two components.
     \item Each hyperplane is itself
      a cubical complex, and inherits a CAT(0) structure
     from $X$.
\end{enumerate}

Let $H$ be a subgroup of the finitely generated group $G$.
The number $e(G,H)$ of {\em ends of $G$ with respect to $H$}
is defined to be the number of ends of the coset graph
of $G$ relative to $H$.

  In \cite{Sageev95}, the following theorem is proven
  
\begin{thm}\label{thm:essential}
Suppose that $G$ acts on a finite dimensional $CAT(0)$ cubical complex 
without a global fixed point. Then there exists a hyperplane $J$ in $X$ such that 
$e(G, stab(J))>1$.
\end{thm}
  
 This result was later generalized
by Niblo and Roller \cite{NibloRoller98} and Gerasimov
\cite{Gerasimov97} to actions on infinite dimensional complexes.

\subsection{The Algebraic Torus Theorem}
We will appeal to the following theorem of Dunwoody and Swenson
\cite{DunwoodySwenson2000}, which is a generalization of a theorem
of Scott and Swarup \cite{ScottSwarup2000} and Bowditch
\cite{Bowditch98b}. The following statement is a slightly weakened
version of the statement appearing in \cite{DunwoodySwenson2000}.

\begin{thm}[The Algebraic Torus Theorem]\label{thm:AAT}
     Suppose that $G$ is a finitely generated group which contains a
     virtually polycyclic subgroup $H$ with $e(G,H)>~1$.
     Then one of the following holds:
    \begin{enumerate}
    \item $G$ is virtually polycyclic.
    \item $G$ has a non-elementary fuchsian quotient
     with virtually polycyclic kernel.
    \item   $G$ splits as a free product with amalgamation or $HNN$-extension  over
     a virtually polycyclic subgroup.
    \end{enumerate}
\end{thm}

\subsection{Groups acting on CAT(0) spaces}
Finally, we record the following important facts
 about groups acting on CAT(0) spaces.

\begin{lem}\label{lem:ss solvable}
(1)
 If $G$ acts cellularly on a CAT(0) complex $X$
with finitely many shapes,
then $G$ acts semi-simply
with a discrete set of translation lengths
\cite{Bridson99}.

\noindent (2)
 If $G$ acts semi-simply and properly-discontinuously
on the CAT(0) space $X$,
then any virtually solvable subgroup $H$ of $G$
is virtually abelian \cite{BridsonHaefliger}.
\end{lem}

Moreover, we have the following theorem (see
\cite[Thm~7.5 and Rem~7.7]{BridsonHaefliger}):
\begin{thm}\label{thm:terminates}
Let $G$ act by semi-simple isometries on a CAT(0) space so
that:
\begin{enumerate}
\item there is a bound on the dimension of an isometrically embedded
flat,
\item the set of translation numbers of elements of $G$  is discrete at
0,
\item there is a bound on the order of finite subgroups.
\end{enumerate}
Then any sequence $H_1\subsetneq H_2\subsetneq \cdots$
of virtually abelain subgroups terminates.
\end{thm}

\begin{cor}\label{cor:conjugate}
Let $G$ act as in Theorem~\ref{thm:terminates}.
If $H^t\subset H$ then $H^t=H$.
\end{cor}
\begin{proof}
This follows from Theorem~\ref{thm:terminates} by setting
$H_i= H^{t^{-i}}$.
\end{proof}

\section{The Main Theorem}
We now prove Theorem~\ref{thm:main} in the following form:
\begin{thm}
     Suppose there is a
     uniform bound on the orders
of finite subgroups of the group $G$.
     Let $G$ act properly on a finite dimensional CAT(0) cubical complex $X$. Then
either $G$ has a
     rank~$2$ free
     subgroup or $G$ is virtually a finitely generated abelian group.
\end{thm}

\begin{proof}
     We proceed by induction on $\dimension(X)$.
     If $\dimension(X)=0$ then $G$ is finite
     and hence virtually finitely generated abelian.
     Suppose the theorem holds for $\dimension(X)<n$,
     and consider the case where $\dimension(X)=n$.

     We first prove
     the theorem in the case that $G$ is finitely generated. If $G$ is finite, we are done, so suppose that $G$ is infinite.
     Since $G$ is infinite and
   acts properly on $X$,
   \begin{com}
   Michah, are properly and properly-discontinously the same for you?
   Could you choose one and adopt uniform language.
   \end{com}
    it has no global fixed point.
     Therefore,
     by Theorem \ref{thm:essential}, there exists a
     hyperplane $J$ in $X$ with $e(G,stab(J))>1$. Let $H=stab(J)$.
     Then $H$ acts properly on $J$. By induction, either
     $H$ is virtually finitely generated abelian or $H$ contains a
rank~$2$
     free subgroup.
     If $H$ has a rank~$2$ free subgroup then so does $G$ and we are
     done. We therefore assume that $H$ is virtually a finitely
generated abelian
     group. Applying Theorem~\ref{thm:AAT} we have either:
            \begin{enumerate}
    \item $G$ is virtually polycyclic,
    \item $G$ has a non-elementary fuchsian quotient with
     virtually polycyclic kernel,
    \item   $G$ splits as a free product with amalgamation or $HNN$-extension  over
     a virtually polycyclic subgroup.
    \end{enumerate}

    In case~(1), $G$ is virtually abelian by Lemma~\ref{lem:ss
    solvable}(2).
    In case (2), the non-elementary fuchsian quotient contains
    a rank 2~free subgroup, and hence $G$ does as well.
    Thus, suppose that
        $G$ splits as
     an amalgamated free product $A*_P B$ or HNN extension $C*_P$
     over a virtually
     polycyclic group $P$.
     Since $X$ has finitely many shapes, $G$ acts semi-simply by
     Lemma~\ref{lem:ss solvable}(1),
     and so $P$ is virtually abelian by Lemma~\ref{lem:ss solvable}(2).

An application of the normal form theorem for graphs of groups
\cite{Serre80}, shows that if either $[A:P]>2$ or $[B:P]>2$ then
$G$ contains a rank~2 free subgroup. We may therefore assume that
$[A:P]=2$ and $[B:P]=2$. In this case $G$ is virtually polycyclic.
Indeed, the Bass-Serre tree is a line, so $G$ has an index~2
subgroup $G'$ which acts by translations, and the resulting
homomorphism $G'\rightarrow \integers$ has kernel~$P$, so $G'\cong
P\semidirect \integers$.

Again, the normal form theorem shows that if both $[C:P]>1$ and
$[C:P^t]>1$ (where $t$ is the stable letter of $C*_P$) then $G$
contains a rank~2 free subgroup. On the other hand, by
Corollary~\ref{cor:conjugate}, if $[C:P]=1$ then $[C:P^t]=1$ (and
vice-versa). We may therefore assume that $[C:P]=1$ and
$[C:P^t]=1$,
  and so $G\cong P\semidirect \integers$.

 In each case, $G$ is virtually polycyclic and hence virtually abelian
 by Lemma~\ref{lem:ss solvable}.

Now suppose that  $G$ is not finitely generated. Thus, $G$ contains an
infinite ascending sequence
of proper subgroups $G_1\subsetneq G_2\subsetneq G_3\dots$.
We are done if any $G_i$ contains a rank~2 free subgroup.
 So we shall show that the assumption that
 each $G_i$ is virtually finitely generated abelian
 leads to a contradiction.
We first verify
conditions~(1)-(3) in the hypothesis of Theorem~\ref{thm:terminates}.
Condition~(1) holds since $X$ is finite dimensional.
Condition~(2) holds by Lemma~\ref{lem:ss solvable},
since $X$ is a finite dimensional cubical complex and
hence a polyhedral complex with finitely many shapes.
Condition~(3) is hypothesized in our theorem.
Thus the infinite sequence $G_1\subsetneq G_2\subsetneq G_3\dots$ terminates by
Theorem~\ref{thm:terminates} which is impossible.
\end{proof}

\bibliographystyle{alpha}

        \bibliography{wise} 


%
%
\end{document}